\documentclass{amsart}
\usepackage{amsmath,amsthm,amssymb,amsfonts}
\usepackage[colorlinks=true]{hyperref}
\usepackage{cleveref}

\hypersetup{urlcolor=blue, citecolor=blue}

\newcommand{\al}{\alpha}    
\newcommand{\de}{\delta}    \newcommand{\De}{\Delta}
  \newcommand{\ep}{\varepsilon}
\newcommand{\la}{\lambda}   
\newcommand{\si}{\sigma}    
\newcommand{\om}{\omega}    
    
\newcommand{\R}{\mathbb{R}}\newcommand{\Z}{\mathbb{Z}}
\newcommand{\N}{\mathbb{N}}

\newcommand{\Sc}{\mathcal{S}}
\newcommand{\F}{\mathcal{F}}\newcommand{\Sp}{\mathbb{S}}
\newcommand{\B}{\mathcal{B}}\newcommand{\Hs}{\dot{\mathcal{H}}}
\def\O{{\mathcal{O}}}

\newcommand{\pt}{\partial_t}\newcommand{\pa}{\partial}
\newcommand{\les}{{\lesssim}}
\newcommand{\beeq}{\begin{equation}}\newcommand{\eneq}{\end{equation}}

\def\<{\langle}             \def\>{\rangle}

\newenvironment{prf}{\noindent {\bf Proof.} }{\endprf\par}
\def \endprf{\hfill  {\vrule height6pt width6pt depth0pt}\medskip}

\newtheorem{thm}{Theorem}[section]
\newtheorem{prp}[thm]{Proposition}
\newtheorem{lem}[thm]{Lemma}

\theoremstyle{remark}
\newtheorem{rem}{Remark}[section]

\theoremstyle{definition}
\newtheorem{defi}{Definition}

\numberwithin{equation}{section}

\begin{document}
\bibliographystyle{plain} 

\title
[nonlinear half wave equation]
{Fractional derivatives of composite functions 
and the Cauchy problem for the nonlinear half wave equation
}

\author{Kunio Hidano}
\address{Department of Mathematics\\
Faculty of Education\\
Mie University\\
1577 Kurima-machiya-cho, Tsu, Mie 514-8507, JAPAN}
\email{hidano@edu.mie-u.ac.jp}

\author{Chengbo Wang}
\address{School of Mathematical Sciences\\                Zhejiang University\\                Hangzhou 310027, China}\email{wangcbo@zju.edu.cn}
\urladdr{http://www.math.zju.edu.cn/wang}

\date{\today}

\begin{abstract}
We show new results of wellposedness for 
the Cauchy problem for the half wave equation 
with power-type nonlinear terms. 
For the purpose, we propose two approaches 
on the basis of the contraction-mapping argument. 
One of them relies upon the $L_t^q L_x^\infty$ Strichartz-type 
estimate together with 
the chain rule of fairly general fractional orders. 
This chain rule has a significance of its own. 
Furthermore, in addition to the weighted fractional chain rule 
established in Hidano, Jiang, Lee, and Wang 
(arXiv:1605.06748v1 [math.AP]), 
the other approach uses weighted space-time $L^2$ estimates 
for the inhomogeneous equation which are recovered 
from those for the second-order wave equation. 
In particular, by the latter approach we settle the problem 
left open in Bellazzini, Georgiev, and Visciglia 
(arXiv:1611.04823v1 [math.AP]) 
concerning the local wellposedness in 
$H^{s}_{{\rm rad}}({\mathbb R}^n)$ with $s>1/2$.
\end{abstract}




\maketitle

\section{Introduction}
This paper is concerned with 
the Cauchy problem for the nonlinear, first-order 
wave equation
\begin{equation}\label{eq-1.1}
\begin{cases}
i\partial_t u-\sqrt{-\Delta}u=F(u),\quad
t>0,\,x\in{\mathbb R}^n,\\
u(0)=u_0.
\end{cases}
\end{equation}
This equation has rich mathematical problems, and 
it has recently gained much attention; 
to mention a few papers treating power-type nonlinear terms, 
see \cite{BGV}, \cite{Dinh}, \cite{FGO} for 
local/global wellposedness, 
\cite{Inui} for finite-time blow up, 
 \cite{KrLeRa13}, \cite{BGV}
 for stability/instability of ground states,
\cite{Dinh} for illposedness for low-regularity data, 
and \cite{FGO} for the proof of various a priori estimates. 
The present paper also treats 
the power-type nonlinear term of the form 
$\lambda|u|^{p-1}u$ (which is algebraic if $p$ odd) or $\lambda|u|^p$  (which is algebraic if $p$ even) 
($\lambda\in{\mathbb C}\setminus\{0\}$, $p>1$), 
and discusses three problems 
left unexplored in the existing literature. 
They are: 
(i) wellposedness in the Sobolev spaces 
$H^{s}(\R^n)$ with 
$s\ge s_c:=n/2-1/(p-1)$, 
especially when $s$ is strictly smaller than, and very close to $p$; 
(ii) local wellposedness in 
$H^{(1/2)+\varepsilon}_{{\rm rad}}({\mathbb R}^n)$ 
when $1<p<1+2/(n-1)$, 
which is the problem left open in \cite{BGV}; 
(iii) global wellposedness for small data 
in $H^s_{{\rm rad}}({\mathbb R}^n)$ 
for $1+2/(n-1)<p<1+2/(n-2)$, $s\in (s_c, 1]$. 

Our proof of wellposedness in 
$H^{s}$ builds upon 
the $L_t^q L_x^\infty$ Strichartz-type estimates 
due to Klainerman and Machedon \cite{KM} 
(see also Proposition 1 of Fang and Wang \cite{FW2}).
Before the statement of the wellposedness,
let us make the assumption on the nonlinear term $F(u)$. 
\begin{defi}\label{def-1}
We say $F(u)$ has the property 
$(F)_{k,p}$ if 
there exists $p>1$ and $k\in{\mathbb N}$ with $k\leq p$ 
such that when it is considered as a function 
${\mathbb R}^2\to{\mathbb R}^2$, 
$F:{\mathbb C}\to{\mathbb C}$ is a $C^k$ function, 
satisfying $F^{(j)}(0)=0$ for $0\leq j\leq k$ and 
for all $z_1,z_2\in{\mathbb C}$
\beeq\label{eq-Leib-assu} |F^{(j)}(z_{1})-F^{(j)}(z_{2})|\le C_{j} 
\left\{
\begin{array}{ll}
|z_{1}-z_{2}|[\max(|z_{1}|, |z_{2}|)]^{p-j-1}      &    j\le k, p\ge j+1 \\
|z_{1}-z_{2}|^{p-j}      &   j=k, p\in [k, k+1]\ .
\end{array}
\right.
\eneq
\end{defi}
Let $U(t):=\exp\{-itD\}$ with $D=\sqrt{-\Delta}$. 
By $\lceil p\rceil$, we mean the smallest integer greater than or equal to $p$.
Our first theorem concernig the problem (i) is as follows:
\begin{thm}\label{th-1.1}
Let $n\ge 1$, 
$p>1$.
When $F(u)$ is not algebraic, 
suppose that 
$$p>\max\left(s_c, \frac{n-1}2,\frac{n+1}4\right)\ ,$$
and
$F(u)$ has the property $(F)_{k,p}$ with
$k=\lceil p\rceil-1$ $($i.e., $k<p\le k+1$ and $k\in \N )$.
Then the problem \eqref{eq-1.1} is locally wellposed in $H^s$, provided that
$s\ge s_c=n/2-1/(p-1)$, and $
s>\max\left((n-1)/2,(n+1)/4\right)$ (we additionally assume $s<p$ when $F$ is not algebraic).
More precisely, we have
\begin{enumerate}
  \item for any $u_0\in H^{s}$, there exists a 
$T_0>0$ depending on $n$, $p$, $u_0$, the constants $C_j$ in \eqref{eq-Leib-assu}, such that there exists a unique solution 
\begin{equation}\label{eq-1.4}
u\in
L^\infty(0,T_0;H^{s})
\cap C([0,T_0];L^2)
\cap 
L^{q}(0,T_0;L^\infty)
\end{equation}
to the associated integral equation
\begin{equation}
u(t)
=
U(t)u_0
-i
\int_0^t
U(t-\tau)
F(u(\tau))d\tau
\end{equation}
 for some $q\in (\max({4}/(n-1), 2),\infty]$ with $q\ge p-1$, $1/q\ge n/2-s$.
 Moreover, $u$ actually belongs to $C([0, T_0]; H^s)$.
  \item {\rm (Persistence of regularity) }
  if $u_0\in H^{s_1}$ with some $s_1>s$ ($s_1< p$ if $F(u)$ is not algebraic),
then $u\in C([0,T_0];H^{s_1})$.
  \item {\rm (Critical wellposedness)} In additon, when we can take $s=s_c$, that is, when
$p>5$ for $n=2$, 
$p>3$ for $n\ge 3$, and if $n\ge 8$ and $F(u)$ is not algebraic, 
$$
p>\frac{n+2+\sqrt{(n+2)(n-6)}}{4}\,\,\,\mathrm{for\,\,}n\geq 8,
$$
so that we have
$p>s_c$, 
the problem is critically local wellposed in $H^{s_c}$. In particular, there exists $\ep_0>0$ such that
$T_0=\infty$ and the solution scatters in $H^{s_c}$, i.e.,
$$\exists u_0^+\in H^{s_c}, \lim_{t\to +\infty}\|u(t)-U(t)u_0^+\|_{H^{s_c}}=0,$$
 when
$\|D^{s_c}u_0\|_{L^2(\R^n)}<\ep_0$. 

\end{enumerate}

\end{thm}

\begin{rem}\label{rem-1.1}
The similar result with additional assumption $p>\lceil s\rceil$ has been obtained in Dinh \cite{Dinh}.
If the initial data $u_0$ is radially symmetric about the origin $x=0$, 
the result similar to  \Cref{th-1.1} obviously remains true for
$s\ge \max(s_c, (n-1)/2)$ for $n\ge 3$
and $s\ge s_c$, $s>1/2$ for $n= 2$.
This improvement is due to the well-known fact that 
the range of admissible pairs $(q,r)$ 
improves for radially symmetric data, see \Cref{th-Stric} below. 
Such a  result of small data global existence for radially symmetric data in the case $n=p=3$ has been already shown in   Fujiwara, Georgiev and  Ozawa \cite{FGO}. 
\end{rem}
\begin{rem}\label{rem-1.2}
The regularity assumptions $s\ge s_c$ and $s>(n+1)/4$ are sharp in general. Actually, by the connection of this problem with the nonlinear wave equations established in \Cref{sec-bu}, we know that,
when $p\ge 2$ with $n=3,4$ and $p>1$ with $n\ge 5$,
there exists $F(u)$ (e.g., $F=i |\Re u|^p$) such that the problem is ill-posed in $H^s$ for $s=\max(s_c, (n+1)/4)-\de$ with arbitrary small $\de>0$. See Lindblad
\cite{Ld93} for $n=3$ and $p=2$, 
Fang and Wang \cite[Theorem 1.2]{FW1} for 
 $s\in (-n/2, s_c)$ when $s_c>0$,
  \cite[Theorem 1.3]{FW1} for
$s\in (\max(s_c, 0), (n+1)/4)$
when $s_c<(n+1)/4$ and $n\ge 5$,
and
  \cite[Corollary 1.1 and Theorem 1.4]{FW1} for
$s\in (s_c, (n+1)/4)$
when $s_c<(n+1)/4$, $p\ge 2$ and $n=3, 4$). 
See also \cite{Dinh} for some ill-posed results with $s<s_c$, by applying the technique of Christ, Colliander and Tao \cite{CCT03}.
\end{rem}

In addition to the $L^q_tL^\infty_x$ Strichartz-type estimate, 
the proof of  \Cref{th-1.1} uses the Ginibre-Ozawa-Velo type 
estimate on fractional derivatives of composite functions 
(see \Cref{thm-Leib-h} below), by which we in particular obtain for some constant 
$C=C(n,p,s)>0$
\begin{equation}\label{eq-1.5}
\|D^s(|v|^{p-1}v)\|_{L^2({\mathbb R}^n)}
\leq
C
\|v\|_{L^\infty({\mathbb R}^n)}^{p-1}
\|D^s v\|_{L^2({\mathbb R}^n)}
\end{equation}
whenever $s\in (0, \min(p, n/2))$, $p>1$. 
The way of showing the fractional chain rule \eqref{eq-1.5} is 
inspired by Ginibre, Ozawa and Velo \cite{GOV}, 
and we closely follow the argument of the proof of 
Lemma 3.4 of \cite{GOV}. 
In the contraction-mapping argument, 
we use  \eqref{eq-1.5} in combination with the basic estimate \Cref{th-Stric} for 
the inhomogeneous equation \eqref{eq-2.2}, 
and hence the $L_t^{q} L_x^\infty$ Strichartz-type norm 
naturally comes into play. 

As mentioned above (see the problem (ii)), 
the second purpose of this paper is to 
prove the local wellposedness in 
$H^{(1/2)+\varepsilon}_{{\rm rad}}({\mathbb R}^n)$, 
when $p$ is $H^{1/2}$-subcritical, that is, $1<p<1+2/(n-1)$. 
The importance of studying the wellposedness in $H^{1/2}$ 
is obvious in view of the two conservation laws 
\begin{align}
&
\|u(t)\|_{L^2({\mathbb R}^n)}
=
\|u_0\|_{L^2({\mathbb R}^n)},\\
&
\frac12\|D^{1/2}u(t)\|_{L^2({\mathbb R}^n)}^2
+
\frac{\lambda}{p+1}
\|u(t)\|_{L^{p+1}({\mathbb R}^n)}^{p+1}\\
&=
\frac12\|D^{1/2}u_0\|_{L^2({\mathbb R}^n)}^2
+
\frac{\lambda}{p+1}
\|u_0\|_{L^{p+1}({\mathbb R}^n)}^{p+1}\nonumber
\end{align}
for $F(u)=\lambda|u|^{p-1}u$ $(\lambda\in{\mathbb R})$. 
In spite of its importance, 
the local wellposedness in $H^{1/2}$ actually remains open, except the particular 
case $n=1$ and $F(u)=-|u|^2 u$,
see 
Krieger, Lenzmann and Rapha\"el 
\cite{KrLeRa13}. 
In the present paper, at the cost of imposing 
a little more regularity and radial symmetry on the data, 
we study the wellposedness in 
$H^{(1/2)+\varepsilon}_{{\rm rad}}({\mathbb R}^n)$. 
We shall prove
\begin{thm}\label{th-1.2}Let $n\geq 2$, $p\in (1,1+2/(n-1))$ and suppose that 
$F(u)$ satisfies 
\begin{equation}\label{eq-1.8}
|F(u)|\leq 
C_0|u|^{p},
|F'(u)|\leq 
C_1 |u|^{p-1}.
\end{equation}
Let $s\in(1/2,1]\cap(1/2,n/2)$ and $s_1\in (1/2, s]$. 
Then there exists a constant $c>0$ depending on 
$n$, $p$, $s$, $s_1$, and the constant $C_j$ in \eqref{eq-1.8}
such that the Cauchy problem  \eqref{eq-1.1}
with $u_0\in H^{s}_{{\rm rad}}({\mathbb R}^n)$ 
admits a unique, radially symmetric solution 
satisfying
\beeq
u\in L^\infty(0,T_\Lambda;H^s)\cap 
C([0,T_\Lambda];L^2),\eneq
\beeq r^{-(1-\delta)/2} D^\sigma u
\in
L^2((0,T_\Lambda)\times{\mathbb R}^n),
\delta:=1-\frac{n-1}{2}(p-1),
\sigma=0,1-s_1, s_1, s,
\eneq
where 
\begin{equation}\label{eq-life-lower}
T_\Lambda:=c\Lambda^{-(p-1)/\delta}, \quad
\Lambda
:=
\|D^{s_1}
u_0
\|_{L^2({\mathbb R}^n)}^{1/2}
\|
D^{1-s_1}
u_0
\|_{L^2({\mathbb R}^n)}^{1/2}
.
\end{equation}
\end{thm}

\begin{rem}  By (1.11), we see that if $\|D^{s_1}u_0\|_{L^2({\mathbb R}^n)}^{1/2}\|D^{1-s_1}u_0\|_{L^2({\mathbb R}^n)}^{1/2}\to 0$  for {\it some} $s_1\in(1/2,s]$, 
then 
$T_\Lambda\to \infty$, 
regardless of $\|u_0\|_{H^s}$. 
\end{rem}

This theorem resolves the problem left open in Bellazzini, Georgiev and Visciglia \cite{BGV}. 
In order to show the wellposedness in 
$H^{(1/2)+\varepsilon}_{{\rm rad}}({\mathbb R}^n)$ 
along with the lower bound (1.11) on lifespan, 
we use the method of weighted space-time $L^2$ estimates 
together with the weighted fractional chain rules 
(see \cite[Theoem 2.5]{HJLW}). 
The similar approach has already played a key role 
in getting the optimal lower bound on lifespan of 
small solutions to 
the second order semi-linear wave equations 
with the minimal-regularity radial data \cite{HJLW}. 
See also Keel, Smith and Sogge \cite{KSS} for an earlier and influential result, 
proved by the approach based on weighted $L^2$-estimates, 
concerning long-time existence of small, classical solutions 
to the second-order semi-linear wave equation.
With the help of the boundedness of 
the Riesz transforms on 
$L^2(\omega(x)dx)$ $(\omega(x)\in A_2)$, 
we recover the space-time $L^2$ estimates 
for the half-wave equation  
from those for the second order wave equation, see \Cref{th-LE} below.

Concerning the problem (iii), we show:
\begin{thm}\label{th-1.3}
Let $n\ge 2$ and let $F(u)$ satisfy  \eqref{eq-1.8}
for some $p\in (1+2/(n-1),1+2/(n-2))$, $s\in (s_c,1]\cap (s_c, n/2)$ and $u_0\in H_{\rm rad}^s({\mathbb R}^n)$. 
Also, let $\delta\in (0,1)$ satisfy 
\begin{equation}
1-\delta
=
\left(
\frac{n}{2}-s
\right)
(p-1)
\end{equation}
and let $\delta'$ satisfy 
\begin{equation}
\delta
<
\delta'
<
\left(
s-\frac12
\right)
(p-1).
\end{equation}
Moreover, set 
\begin{equation}
s_1
:=
\frac{n}{2}-\frac{1-\delta+\delta'}{p-1}
\end{equation}
so that $s_1\in (1/2,s_c)$. 
Then there exists an $\varepsilon_0>0$ such that 
if 
\begin{equation}
\|
D^{s_1}
u_0
\|_{L^2({\mathbb R}^n)}
+
\|D^s u_0\|_{L^2({\mathbb R}^n)}
<
\varepsilon_0,
\end{equation} 
then 
the Cauchy problem \eqref{eq-1.1}  admits a unique, radially symmetric solution 
satisfying 
\begin{align}
&
u\in L^\infty(0,\infty;H^s)\cap 
C([0,\infty);L^2),\\
&
r^{-(1-\delta)/2}
\langle r\rangle^{-\delta'/2}
D^\sigma u
\in
L^2((0,\infty)\times{\mathbb R}^n),\,\,
\sigma=0,s_1,s.
\end{align}
\end{thm}

\begin{rem} As mentioned in \Cref{rem-1.1},
global existence of small solutions for radially symmetric data in $H^{s_c}_{\mathrm{rad}}$ can be 
proved by the method using the Strichartz estimates 
for $n=2$, $p>3$.
\end{rem}

Our proof of 
\Cref{th-1.3}
 uses 
global-in-time space-time $L^2$ estimates 
with the weight 
$r^{-(1-\delta)/2}\langle r\rangle^{-\delta'/2}$, 
as in the study of the second order wave equations \cite{HWY2}, 
\cite{HJLW}. See also the earlier papers \cite{HiYo06}, \cite{HWY1} 
where the weight 
$r^{-(1-\delta)/2}\langle r\rangle^{-\delta/2}$ 
was used for the proof of ``almost global'' existence.

We note that 
in view of \Cref{th-1.1} and \Cref{th-1.3}, 
global existence of small solutions 
remains open when $1+2/(n-2)\leq p<3$ $(n\geq 4)$. We expect that 
one may improve the result to the cases $p>s_c$ and $p>1+2/(n-1)$,
with further efforts. However, for the cases with $p<s_c$, e.g., $p\in (7/5, 3/2)$ with $n=7$, the problem seems to be beyond the scope of the current technology.

We should also mention that 
Inui  \cite[Theorem 1.4]{Inui} proved finite-time blow up 
even for small data in the case $F(u)=\lambda |u|^p$ with $1<p<1+2/n$. 
Since we prove global existence 
for $1+2/(n-1)<p<1+2/(n-2)$ 
and small $H^1$ radially symmetric data, 
it is an interesting problem 
whether \eqref{eq-1.1} with certain $F(u)$ with power
$1+2/n\leq p\leq 1+2/(n-1)$ has finite-time blow-up solutions even for 
small data. 
In this connection, 
it is also an interesting problem whether the lower bound  of the lifespan \eqref{eq-life-lower} is sharp. 
 By establishing a connection between the half-wave problems and the nonlinear wave equations, we show that the our lower bound is optimal, in general. Moreover, we could also handle the critical case.
\begin{thm}\label{th-1.4} Let  $n=2, 3, 4$, 
 $F(u)=i |\Re u|^p$ with $p\in (\max(1, (n-1)/2),1+2/(n-1)]$, and
$s\in [(n-1)/2, p)$ except $s=1/2$ and $n=2$.
For any real-valued, radial function $g\in C_0^\infty(\R^n)$ with $\int gdx>0$, 
there exist constants $C, \ep_0>0$ such that
for any $\ep\in (0, \ep_0)$, the Cauchy problem \eqref{eq-1.1}
with data $u_0=\ep g$ can not have solution $u\in
L^\infty(0,L_{\ep};H^s)\cap 
C([0,L_{\ep}];L^2)$
with
\beeq\label{eq-upperbd}
L_{\ep}= \left\{
\begin{array}{ll} 
  \exp(C \ep^{-(p-1)}) & p=1+\frac{2}{n-1}\\
  C\ep^{-\frac{2(p-1)}{2-(n-1)(p-1)}} & 1<p<1+\frac{2}{n-1}\ .\end{array}
\right.\eneq
In addition, let $n\ge 2$, $s\in(1/2,1]\cap(1/2,n/2)$ and let $F(u)$ satisfy \eqref{eq-1.8} with $p=1+2/(n-1)$.
There exist $c, \ep_1>0$, such that for any
$\ep\in (0, \ep_1)$, the Cauchy problem \eqref{eq-1.1} admits a solution in 
$u\in
L^\infty(0,T_\ep;H^s)\cap 
C([0,T_\ep];L^2)$, where
\beeq\label{eq-life-lower2}
T_\ep=
  \exp(c\ep^{-(p-1)}) 
\eneq
for any
radial data
$u_0 \in H^s$ with 
$
\ep^2=
  \| u_0\|^{2}_{\dot H^{s}}+ \|u_0\|_{\dot H^{s}}
\|u_0\|_{\dot H^{1-s}} 
$.
\end{thm}

This paper is organized as follows.
Section \ref{sec-chain}  is devoted to the proof of a
fractional chain rule for general composite functions,
 \Cref{thm-Leib-h}. 
 In Section \ref{sec-st}, 
for half wave equations,
we collect and study two class of space-time estimates, Strichartz type estimates, \Cref{th-Stric}, and Morawetz type local energy estimates, \Cref{th-LE}.
Such estimates have been extensively investigated for the wave equations, and have played an important role in our understanding of various linear and nonlinear wave equations.
Equipped with the Strichartz type estimates,
 and 
the chain rule of fractional orders,
we present the proof of
\Cref{th-1.1} in \Cref{sec-1.1}.
In \Cref{sec-1.2},
we prove Theorems \ref{th-1.2} and \ref{th-1.3}, as well as the existence part in Theorem \ref{th-1.4}.
The key to the proof is the local energy estimates \Cref{th-LE}, the radial Sobolev inequalities, and a weighted fractional chain rule of  \cite[Theorem 2.5]{HJLW}.
\Cref{sec-bu} is devoted to the remaining part of  Theorem \ref{th-1.4},
nonexistence of global solutions and upper bound of the lifespan, for the sample case of $F(u)=i |\Re u|^p$ with $1<p\le 1+2/(n-1)$.

\section{Fractional chain rule}\label{sec-chain}
In this section, we prove a version of the fractional chain rule, which has a significance of its own and
plays an essential role in the proof of Theorem \ref{th-1.1}. 

\subsection{Function spaces}
To begin with, let us recall some basic facts about the Besov/Sobolev spaces.
At first,
let $\phi(x)\in C_0^\infty(\R^n)$, $\phi\ge 0$, with $\phi=1$ for $|x|\le 1$ and $\phi=0$ for $|x|\ge 2$, we define the Littlewood-Paley projection operators $S_j$, $P_j$, for $f\in \Sc'$ as follows
$$\F (S_j f)(\xi)=\phi(2^{-j} \xi)(\F f)(\xi)\ , P_j f=S_{j}f-S_{j-1}f\ ,$$
where $\F$ denotes the Fourier transform.
We also define the inhomogeneous Littlewood-Paley projection operators by $\tilde P_j$ for $j\in\N$ as follows $\tilde P_j=P_j$ for $j\ge 1$ and  $\tilde P_0=S_0$. It is well known that for any $f\in \Sc'$, we have
$$\lim_{N\rightarrow +\infty}S_N f=\lim_{N\rightarrow +\infty} \sum_{0\le j\le N}\tilde P_j f= f\ \mathrm{in}\ \Sc'$$
However, it is not always true for the homogeneous operators. For that purpose, we introduce (see \cite[Definition 1.26, page 22]{BCD})
$$\Sc_h':=\{f\in \Sc'; \lim_{N\rightarrow -\infty}\|S_N f\|_{L^\infty}=0\}\ .$$
It is obvious that $\Sc_h'$ is included in the space of tempered distributions vanishing at infinity,
$$\Sc_0':=\{f\in \Sc'; \lim_{N\rightarrow -\infty}S_N f=0\ \mathrm{in}\ \Sc'\}\ ,$$
where we have $f=\sum_{j\in \Z} P_j f$ in $\Sc'$.

For $s\in \R$, $q,r\in [1,\infty]$, we define the Besov (semi)norms for $f\in \Sc'$,
\beeq\label{def-Besov-1}\|f\|_{\dot B^{s}_{q,r}}=\|2^{js}P_{j}f\|_{\ell^{r}_{j\in\Z}L^{q}}\ ,
\|f\|_{B^{s}_{q,r}}=\|2^{js}\tilde P_{j}f\|_{\ell^{r}_{j\ge 0}L^{q}}\ .
\eneq Based on the Besov (semi)norm, we define Besov spaces by
\beeq
B^{s}_{q,r}=\{f\in\Sc';
\|f\|_{B^{s}_{q,r}}<\infty
\}\ ,\ 
\dot \B^{s}_{q,r}=\{f\in\Sc'_h;
\|f\|_{\dot B^{s}_{q,r}}<\infty
\}\ .\eneq
It is known that the  inhomogeneous Besov spaces are Banach spaces, for all $s\in \R$ and $q,r\in [1,\infty]$, and the homogeneous Besov spaces $\dot\B^s_{q,r}$ are Banach spaces,  if and only if
\beeq \label{Besov-complete}s<n/q, r\in [1,\infty] \textrm{ or } s=n/q, r=1\ ,\eneq
see, e.g., \cite[Theorem 2.25 and Remark 2.26, pages 67-68]{BCD}.
Moreover, if $\|f\|_{\dot \B^{s}_{q,r}}<\infty$ for $(s,q,r)$ with \eqref{Besov-complete} satisfied, the series
$\sum_j P_j f$ is convergent to some $g\in\Sc_h'$, see, e.g., \cite[Remark 2.24, page 66]{BCD}.

Similarly, for $s\in \R$, we define the  $L^2$ based Sobolev (semi)norms for $f\in \Sc'$,
\beeq\label{def-Sob-1}\|f\|_{\dot H^{s}}=\||\xi|^s \F(f)(\xi)\|_{L^{2}}\ ,
\|f\|_{H^{s}}=\|(1+|\xi|^2)^{s/2} \F(f)(\xi)\|_{L^{2}}\ .
\eneq Based on these (semi)norms, we define $L^2$ based Sobolev spaces by
\beeq
H^{s}=\{f\in\Sc', 
\|f\|_{H^{s}}<\infty
\}\ ,\ 
\Hs^{s}=\{f\in\Sc', \F(f)\in L^1_{loc}, \|f\|_{\dot H^{s}}<\infty\}\ .\eneq
It is known that the  inhomogeneous Sobolev spaces are Hilbert spaces, for all $s\in \R$, and the homogeneous Sobolev spaces $\Hs^s$ are Hilbert spaces,  if and only if
$s<n/2$,
see, e.g., \cite[Proposition 1.34, page 26]{BCD}.
Moreover, $\dot\B^{s}_{2,2}=\Hs^s$ if $s<n/2$ and we denote
$\|f\|_{\dot\B^{s}_{q,r}}=\|f\|_{\dot B^{s}_{q,r}}$ with property $f\in \dot\B^{s}_{q,r}$.

Here, for future reference, let us record the fractional Leibniz rule for such spaces, see, e.g., \cite[Corollary 2.54]{BCD}.
\begin{lem}[Fractional Leibniz rule]\label{thm-leibniz}
Let $s>0$ with \eqref{Besov-complete}.
 Then $\dot \B^{s}_{q,r}\cap L^\infty$ is an algebra.
 Moreover, there exists a constant $C$, depending only on the dimension $n$, such that
$$\|fg\|_{\dot \B^{s}_{q,r}}\le \frac{C^{s+1}}s\left( \|f\|_{L^{\infty}} \|g\|_{\dot \B^{s}_{q,r}}+\|f\|_{\dot \B^{s}_{q,r}}\|g\|_{L^{\infty}}\right)\ .$$
In particular, we have \beeq\label{eq-Leib}\|fg\|_{\Hs^{s}}\les \|f\|_{L^{\infty}} \|g\|_{\Hs^{s}}+\|f\|_{\Hs^{s}}\|g\|_{L^{\infty}}\ , \ 0<s<n/2,\eneq
and for algebraic $F(u)$ with power $p>1$, 
\beeq\label{eq-Leib2}
    \|F( u)\|_{\Hs^{s}}\les
    \|u\|_{L^{\infty}}^{p-1}    \|u\|_{\Hs^{s}}
    ,\ 0<s<n/2.
\eneq 
\end{lem}

For Besov spaces, we recall the following difference characterization,
see, e.g., \cite[Theorem 6.2.5, 6.3.1]{BL}. Notice, however, that our statement for the homogeneous Besov space is more precise, as the proof of \eqref{eq-Besov-1'} in \cite[Theorem 6.3.1]{BL} has implicitly used the assumption $f=\sum P_j f$, that is, $f\in\Sc'_0$.
\begin{lem}[Difference characterization of the Besov spaces]\label{thm-lem-diff}
For $s\in (0,m)$,  $m\in\N$, $q,r\in [1,\infty]$, we have for any $u\in \Sc'$,
\beeq\label{eq-Besov-1}\|u\|_{\dot B^{s}_{q,r}}\les
\|t^{-s}\|\Delta_m^y u(x)\|_{L^\infty_{y\in B_t} L^{q}_{x}}\|_{L^{r}((0,\infty), dt/t)}\ ,
\eneq
\beeq\label{eq-Besov-2}\|u\|_{B^{s}_{q,r}}
\simeq 
\|u\|_{L^{q}}+
\|t^{-s}
\|\Delta_m^y u(x)\|_{L^\infty_{y\in B_t} L^{q}_{x}}
\|_{L^{r}((0,\infty), dt/t)}\ ,\eneq
where
\beeq\label{eq-mod-0}
\Delta_m^y u(x)=(\tau_{y}-I)^m u(x)\ ,
(\tau_{y}-I)u(x)=u(x+y)-u(x)\ .
\eneq 
Moreover, if $u\in \Sc'_h$, we have
\beeq\label{eq-Besov-1'}\|u\|_{\dot B^{s}_{q,r}}\gtrsim
\|t^{-s}\|\Delta_m^y u(x)\|_{L^\infty_{y\in B_t} L^{q}_{x}}\|_{L^{r}((0,\infty), dt/t)}\ .
\eneq
\end{lem}

\subsection{Fractional chain rule}
With the definitions of the Besov/Sobolev spaces, we are ready to state the fractional chain rule.
\begin{lem}\label{thm-Leib-h}
Suppose that $F(u)$ has the property $(F)_{k,p}$ $($see \Cref{def-1}$)$. 
 Then for any $s\in (0,\min(k+1,p))$ and $q, r\in [1,\infty]$,
we have
\beeq\label{eq-chain-1}
    \|F( u)\|_{\dot B^{s}_{q,r}}\les
    \|u\|_{L^{\infty}}^{p-1}    \|u\|_{\dot B^{s}_{q,r}}
    ,\  u\in \dot \B^{s}_{q,r}\cap L^\infty, 
\eneq 
\beeq\label{eq-chain-2}
    \|F( u)\|_{B^{s}_{q,r}}\les
    \|u\|_{L^{\infty}}^{p-1}    \|u\|_{B^{s}_{q,r}},
    \  u\in  B^{s}_{q,r}\cap L^\infty\ .
\eneq 
In addition, if $r\in [1,2]$, $q\in [1,\infty)$, $s\in (0, n/q)$ and $s\ge n/q-n/2$, then
$F(u)\in\Sc_h'$ and so $F(u)\in \dot \B^{s}_{q,r}$.
In particular, with $q=r=2$, we have
\beeq\label{eq-chain-3}
    \|F( u)\|_{\Hs^{s}}\les
    \|u\|_{L^{\infty}}^{p-1}    \|u\|_{\Hs^{s}}
    , s\in (0,\min(k+1,p,n/2)),\ 
\eneq 
\beeq\label{eq-chain-4}
    \|F( u)\|_{H^{s}}\les
    \|u\|_{L^{\infty}}^{p-1}    \|u\|_{H^{s}}\ ,
        s\in (0,\min(k+1,p))\ .
\eneq 
\end{lem}
\begin{rem}
The estimates of this form have been well-known for many special situations.
Typical examples
include $0<s\le [p]$ or $k=1$. See \cite[Lemma 3.4]{GOV} for the case $k=1$ (actually, our proof is inspired by the proof in \cite{GOV}).
The estimates for $H^s$ ($s>0$) or $W^{s,q}$ ($s\in \N$) have been well-known for the case $s\le [p]$, see e.g. \cite[Lemma
A.9]{Tao06},
\cite[Corollary 6.4.5]{Ho97}.
Estimates involving $L^p$ based Sobolev spaces are also known for may cases, 
see, e.g., 
 \cite[Lemma A2]{Kato95}, 
\cite[Chapter 2 Proposition 5.1]{Tay00} and references therein for case $0<s<1<p$. See also \cite[Theorem 2.5]{HJLW} for a recent work on weighted fractional chain rule.
There are also estimates for H\"older continuous case $0<s<p<1$, see \cite[Lemma A.12]{KV}.
\end{rem}
\begin{rem}
The inhomogeneous estimates \eqref{eq-chain-2} and
\eqref{eq-chain-4} have been known, see, e.g., \cite[Theorem 2, page 325]{RS96}.
It seems to the authors that the homogeneous estimates \eqref{eq-chain-1} and \eqref{eq-chain-3} may be new, although a weaker version of \eqref{eq-chain-1}, with $L^{\infty}$
replaced by $L^{\infty}\cap \dot B^{0}_{\infty, 2}$ and $q,r\ge 2$, has been available from
\cite[Lemma 2.2, page 402]{NO97}, which was also proven by enhancing the proof in \cite{GOV} for $k=1$.
\end{rem}

\subsection{Preparation}
The basic idea of proof is to exploit the
difference characterization of the Besov spaces,
\Cref{thm-lem-diff}. For that purpose, we introduce some notations and exploit the properties of the difference operators.

For fixed $y\in \R^{n}$ with $|y|<t$, let
$$u_{(m)}(x,\la)=[\Pi_{j=1}^{m}(I+\la_{j}\De_{1}^{y})]u(x)
=\sum_{\al\in \{0,1\}^{m}}\la^{\al} \De_{|\al|}^{y} u(x), \la\in [0,1]^{m}\ .
$$
By definition, we see that
$\tau_y \De_{m}^{y}u(x)=\De_{m}^{y}u(x)+\De_{m+1}^{y}u(x)$, and so
\beeq\label{eq-mod-1'}\sup_{|y|<t} |\tau_y \De_{m}^{y}u(x)|\le
|\De_{m}^{y}u(x)|+|\De_{m+1}^{y}u(x)|
, m\ge 1,\eneq
\beeq\label{eq-mod-2}
\|u_{(m)}(x,\la)\|_{L^\infty_x}\les \|u\|_{L^\infty_x}, \la\in [0,1]^m\ , m\ge 0\ .
\eneq
Moreover, for any  $\la\in[0,1]^{m}$, $|y|<t$,
 $k\ge 1$ and $m\ge 0$,
 we have
\beeq\label{eq-mod-3}\De_{k}^{y} u_{(m)}(x,\la)
=\sum_{\al\in \{0,1\}^{m}}\la^{\al} \De_{k+|\al|}^{y} u(x)=
\O(
\sum_{j=0}^m |\De_{k+j}^{y} u(x)|
)\ ,\eneq
\beeq\label{eq-mod-3'}\tau_y\De_{k}^{y} u_{(m)}(x,\la)
=\sum_{\al\in \{0,1\}^{m}}\la^{\al} \tau_y\De_{k+|\al|}^{y} u(x)=
\O(
\sum_{j=0}^{m+1} |\De_{k+j}^{y} u(x)|
)\ .\eneq
For any $\al\in \{0,1\}^{m}$, let $A_{\al}= \{j: \al_{j}=0\}$,
\beeq\label{eq-mod-4}\pa_{\la}^{\al}u_{(m)}(x,\la)=\Delta_{|\al|}^{y} u_{(m-|\al|)}(x, \{\la_{j}\}_{ j\in A_{\al}})\ .\eneq

Notice that \begin{eqnarray*}
\De_{1}^{y} (F\circ u)(x) & = & F(u(x+y))-F(u(x)) \\
 & = & \int_{0}^{1}
\frac{d}{d\la} F(u(x)+\la \De_{1}^{y}u(x))d\la \\
 & = & \int_{0}^{1}
\frac{d}{d\la} F(u_{(1)}(x,\la))d\la \ ,
\end{eqnarray*} 
then for any $k\ge 1$, we have
\beeq
\De_{k}^{y} (F\circ u)(x)=\int_{[0,1]^{k}}
 \frac{d^{k}}{d\la_{1}d\la_{2}\cdots d\la_{k}} F(u_{(k)}(x,\la))d\la \label{eq-comp-1}\ .
\eneq
By using \eqref{eq-mod-4}, we get
\begin{eqnarray*}
&&\De_{k}^{y} (F\circ u)(x) \\
& = & 
\sum_{\substack{1\le l\le k,  \sum \al_{m}=\{1\}^{k} \\
\al_{m}\in \{0,1\}^{k}\backslash 0}}
\int_{[0,1]^{k}}
F^{(l)}(u_{(k)}(x,\la))\Pi_{m=1}^{l}
\De_{|\al_{m}|}^y u_{{(k-|\al_{m}|)}}(x, \{\la_{j}\}_{ j\in A_{\al_{m}}})
d\la
\\ & = & 
\int_{[0,1]^{k}}
F^{(k)}(u_{(k)}(x,\la))\Pi_{m=1}^{k}
\De_{1}^y u_{{(k-1)}}(x, \{\la_{j}\}_{ j\neq m})
d\la\\
&&+
\sum_{\substack{1\le l< k\\  \sum \al_{m}=\{1\}^{k} \\
\al_{m}\in \{0,1\}^{k}\backslash 0}}
\int_{[0,1]^{k}}
F^{(l)}(u_{(k)}(x,\la))\Pi_{m=1}^{l}
\De_{|\al_{m}|}^y u_{{(k-|\al_{m}|)}}(x, \{\la_{j}\}_{ j\in A_{\al_{m}}})
d\la\ .
\end{eqnarray*}
Based on this fact, we know that
\begin{eqnarray*}
&&\De_{k+1}^{y} (F\circ u)(x) 
\\ & = &
\De_1^y \int_{[0,1]^{k}}
F^{(k)}(u_{(k)}(x,\la))\Pi_{m=1}^{k}
\De_{1}^y u_{{(k-1)}}(x, \{\la_{j}\}_{ j\neq m})
d\la \nonumber \\
&&+ \De_{1}^{y}\sum_{\substack{1\le l< k\\ \sum \al_{m}=\{1\}^{k} \\
\al_{m}\in \{0,1\}^{k}\backslash 0}}
\int_{[0,1]^{k}}
F^{(l)}(u_{(k)}(x,\la))\Pi_{m=1}^{l}
\De_{|\al_{m}|}^y u_{{(k-|\al_{m}|)}}(x, \{\la_{j}\}_{ j\in A_{\al_{m}}})
d\la 
\\ & = &
\int_{[0,1]^{k}}
\left(F^{(k)}(u_{(k)}(x+y,\la))-
F^{(k)}(u_{(k)}(x,\la))\right)\Pi_{m=1}^{k}
\De_{1}^y u_{{(k-1)}}(x, \{\la_{j}\}_{ j\neq m})
d\la
\nonumber   \\
 &  & 
+\sum_{l=1}^k \int_{[0,1]^{k}}
F^{(k)}(u_{(k)}(x+y,\la))
 \left(\Pi_{m=1}^{l-1}
\De_{1}^y u_{{(k-1)}}(x+y, \{\la_{j}\}_{ j\neq m})\right)
\\
&&
\times\De_{2}^y u_{{(k-1)}}(x, \{\la_{j}\}_{ j\neq l})
 \left(\Pi_{m=l+1}^{k}
\De_{1}^y u_{{(k-1)}}(x, \{\la_{j}\}_{ j\neq m})\right)
d\la\\
&&+ \sum_{\substack{1\le l< k\\ \sum \al_{m}=\{1\}^{k} \\
\al_{m}\in \{0,1\}^{k}\backslash 0}}
\int_{[0,1]^{k}}
(\De_{1}^{y}F^{(l)}(u_{(k)}(x,\la)))\Pi_{m=1}^{l}
\De_{|\al_{m}|}^y u_{{(k-|\al_{m}|)}}(x, \{\la_{j}\}_{ j\in A_{\al_{m}}})
d\la \\
&&+ \sum_{\substack{1\le l< k\\ \sum \al_{m}=\{1\}^{k} \\
\al_{m}\in \{0,1\}^{k}\backslash 0}}
\sum_{l_0=1}^l
\int_{[0,1]^{k}}
F^{(l)}(u_{(k)}(x+y,\la)) \left(\Pi_{m=1}^{l_0-1}
\De_{|\al_{m}|}^y u_{{(k-|\al_{m}|)}}(x+y, \{\la_{j}\}_{ j\in A_{\al_{m}}})\right)
\\
&&\times
\De_{|\al_{l_0}|+1}^y u_{{(k-|\al_{l_0}|)}}(x, \{\la_{j}\}_{ j\in A_{\al_{l_0}}})
 \left(\Pi_{m=l_0+1}^{l}
\De_{|\al_{m}|}^y u_{{(k-|\al_{m}|)}}(x, \{\la_{j}\}_{ j\in A_{\al_{m}}})\right)
d\la \\
 & = &
\int_{[0,1]^{k}}
\left(F^{(k)}(u_{(k)}(x+y,\la))-
F^{(k)}(u_{(k)}(x,\la))\right)\Pi_{m=1}^{k}
\De_{1}^y u_{{(k-1)}}(x, \{\la_{j}\}_{ j\neq m})
d\la
\nonumber  \\
&&+ \O\left(\sum_{1\le l\le k,  \sum a_{m}\ge k+1 , a_m\in [1, k+1]}
\|u\|_{L^\infty}^{p-l}\Pi_{m=1}^{l}
|\De_{a_{m}}^y u( x)|
\right)\ ,
\end{eqnarray*}
where we have used the assumption \eqref{eq-Leib-assu} with $z_2=0$, \eqref{eq-mod-1'}-\eqref{eq-mod-3'}.

\subsection{Proof of Lemma \ref{thm-Leib-h}}

With help of the previous properties, we could apply Lemma \ref{thm-lem-diff} to give the proof of Lemma \ref{thm-Leib-h}. Here we write down the proof for \eqref{eq-chain-1} only and leave the similar proof for \eqref{eq-chain-2} to the interested reader.

Actually, 
by  \eqref{eq-Leib-assu} and
\eqref{eq-mod-2}-\eqref{eq-mod-3}, we have 
\begin{eqnarray*}&&
\left|F^{(k)}(u_{(k)}(x+y,\la))-
F^{(k)}(u_{(k)}(x,\la))\Pi_{m=1}^{k}
\De_{1}^y u_{{(k-1)}}(x, \{\la_{j}\}_{ j\neq m})\right|
\\
 & \le & 
C(\sum_{j=1}^k |
\De_j^y
u(x)|)^k
\left\{
\begin{array}{ll}
\|u\|_{L^\infty}^{p-k-1}    |\De_1^y u_{(k)}(x,\la)|
  &   p\ge k+1 \\
|\De_1^y u_{(k)}(x,\la)|^{p-k}      &   k\le p\le k+1
\end{array}
\right. 
  \\
 & \le & 
C
\left\{
\begin{array}{ll}
\|u\|_{L^\infty}^{p-k-1}  
(\sum_{j=1}^{k+1} |\De_{j}^y u(x)|)^{k+1}
      &   p\ge k+1 \\
(\sum_{j=1}^{k+1} |\De_{j}^y u(x)|)^{p}     &   k\le p\le k+1
\end{array}
\right. 
\end{eqnarray*}
for any $\la\in [0,1]^m$.
In combination with the  previous relation, we get
\beeq\label{eq-Fu-diff1}\De_{k+1}^y(F(u))(x)\les
\sum_{\substack{1\le l\le k+1\\  \sum a_{m}\ge k+1\\ a_m\in [1, k+1]}}
\|u\|_{L^\infty}^{p-l}\Pi_{m=1}^{l}
|\De_{a_{m}}^y u( x)|, p\ge k+1,
\eneq and
\beeq\label{eq-Fu-diff2}\De_{k+1}^y(F(u))(x)\les
\sum_{j=1}^{k+1} |\De_{j}^y u(x)|^{p} +
 \sum_{
\substack{1\le l\le k, \\ \sum a_{m}\ge k+1 ,\\ a_m\in [1, k+1]}}
\|u\|_{L^\infty}^{p-l}\Pi_{m=1}^{l}
|\De_{a_{m}}^y u(x)|
\eneq when
$k\le p\le k+1$.


For the first case, $0<s<k+1\le p$, we observe that
for any $\sum_{m=1}^{l} a_{m}\ge k+1$ and $a_m\in [1, k+1]$ with $l\in [1, k+1]$,
there exists $b_m\in [1, a_m]$ such that
$\sum_{m=1}^{l} b_{m}= k+1$ and
$0<s b_m/(k+1)<b_m\le a_m$.
Using the difference characterization of the Besov space, \eqref{eq-Besov-1}, we get
\begin{eqnarray*}
&&\|F(u)\|_{\dot B^s_{q,r}} \\& \les & 
\|t^{-s}\|\De_{k+1}^y(F(u))(x)\|_{L^\infty_{y\in B_t}L^{q}_{x}}\|_{L^{r}((0,\infty), dt/t)}
 \\
& \les & 
\sum_{\substack{1\le l\le k+1, \\ \sum a_{m}\ge k+1 ,\\ a_m\in [1, k+1]}}
\|u\|_{L^\infty}^{p-l}
\|t^{-s} \|
\Pi_{m=1}^{l} |\De_{a_{m}}^y u( x)|
\|_{L^\infty_{y\in B_t} L^{q}_{x}}\|_{L^{r}((0,\infty), dt/t)}
\\
& \les & 
\sum_{\substack{1\le l\le k+1, \\ \sum a_{m}\ge k+1 ,\\ a_m\in [1, k+1]}}
\|u\|_{L^\infty}^{p-l}
\Pi_{m=1}^{l} \|t^{-\frac{s b_m}{k+1}} \|
\De_{a_{m}}^y u( x)
\|_{L^\infty_{y\in B_t}L^{\frac{q (k+1)}{b_m}}_{x}}\|_{L^{\frac{r(k+1)}{b_m}}((0,\infty), dt/t)}
\\
& \les & 
 \sum_{\sum_{m=1}^{l} b_{m}=k+1 , b_m\ge 1}
\|u\|_{L^\infty}^{p-l}
\Pi_{m=1}^{l} \|u\|_{\dot B^{\frac{s b_m}{k+1}}_{\frac{q (k+1)}{b_m}, \frac{r(k+1)}{b_m}}}
\\
& \les & 
 \sum_{\sum_{m=1}^{l} b_{m}=k+1 , b_m\ge 1}
\|u\|_{L^\infty}^{p-l}
\Pi_{m=1}^{l} \|u\|_{\dot B^{s}_{q , r}}^{ b_m/(k+1)}\|u\|_{\dot B^{0}_{\infty, \infty}}^{1-b_m/(k+1)}
\\
& \les & 
\|u\|_{L^\infty}^{p-1}
 \|u\|_{\dot B^{s}_{q , r}}\ ,\end{eqnarray*}
 where we have used the interpolation inequality
$$\|u\|_{\dot B^{s\theta}_{q/\theta , r/\theta}}\les \|u\|_{\dot B^{s}_{q, r}}^\theta \|u\|_{\dot B^{0}_{\infty , \infty}}^{1-\theta}, \theta\in (0,1), s>0, q,r\in [1,\infty]\ ,$$ and the trivial embedding
$\|u\|_{\dot B^{0}_{\infty, \infty}}\les \|u\|_{L^\infty}$.

For the second case, $0<s<p$ and $p\in [k, k+1]$, 
since
$0<s<p\le k+1$, all terms except
$\sum_{j=1}^{k+1} |\De_{j}^y u(x)|^{p}$
could be handled in the same way.  Then, we have
\begin{eqnarray*}
\|F(u)\|_{\dot B^s_{q,r}} & \les & 
\|u\|_{L^\infty}^{p-1}
 \|u\|_{\dot B^{s}_{q , r}}+
\sum_{j=1}^{k+1}\|t^{-s}\|
|\De_{j}^y u(x)|^{p}
\|_{L^\infty_{y\in B_t} L^{q}_{x}}\|_{L^{r}((0,\infty), dt/t)}
 \\
& \les & 
\|u\|_{L^\infty}^{p-1}
 \|u\|_{\dot B^{s}_{q , r}}+
\sum_{j=1}^{k+1}\|t^{-s/p}\|
\De_{j}^y u(x)
\|_{L^{qp}_{x}}\|_{L^{rp}((0,\infty), dt/t)}^p
\\
& \les & 
\|u\|_{L^\infty}^{p-1}
 \|u\|_{\dot B^{s}_{q , r}}+
\|u\|_{\dot B^{s/p}_{qp,rp}}^p
\\
& \les & 
\|u\|_{L^\infty}^{p-1}
 \|u\|_{\dot B^{s}_{q , r}}\ .\end{eqnarray*}
 
To complete the proof of \Cref{thm-Leib-h}, it remains to prove that
$F(u)\in\Sc_h'$ if $r\in [1,2]$, $q\in [1,\infty)$, $s\in (0, n/q)$ and $s\ge n/q-n/2$.
Actually,  under the additional assumption on $(s,q,r)$, as $\|F(u)\|_{\dot B^s_{q,r}} <\infty$,
we see that
there is $w\in \Sc_h'$ such that
$$w=\sum_j P_j F(u)\in 
 \dot \B^{s}_{q,r}\subset
\Sc_h' $$
(see, e.g., \cite[Remark 2.24, page 66]{BCD}), which gives us 
$F(u)-w=P(x)$ for some polynomial $P(x)$.
By Sobolev embedding, 
for $q_0\in [2,\infty)$ with $\frac{n}{q_0}=\frac{n}{q}-s$,
we get
$$u, w\in \dot \B^{s}_{q,r}\subset  \dot\B^{0}_{q_0,r}\subset  \dot\B^{0}_{q_0,2}\subset L^{q_0}\ .$$ 
As $F(0)=0$, $u\in L^\infty$ and $F'(u)\in L^\infty$, then $F(u)\in L^{q_0}$. Thus
$$P(x)=F(u)-w\in L^{q_0}$$
and so $P(x)\equiv 0$, which proves $F(u)=w\in\Sc_h'$.
 
\section{Space-time estimates for half wave equations}\label{sec-st}
In this section, 
for half wave equations:
\begin{equation}\label{eq-2.2}
i\partial_t u-\sqrt{-\Delta}u=G,
t>0,\, u(0)=u_0\,\ ,
\end{equation}
we collect and study two class of space-time estimates, Strichartz type estimates and Morawetz type estimates. For the wave equations, these estimates have been extensively investigated and have played an important role in our understanding of various linear and nonlinear wave equations.

The Strichartz type estimates for the half wave equations 
could be easily covered by what have been developed for wave equations. On the other hand, it is less trivial that many Morawetz type estimates could also be recovered from the corresponding estimates for wave equations.

\subsection{Strichartz type estimates} 

We need the $L_t^qL_x^\infty$ Strichartz type estimates. 
More precisely, we use:
\begin{lem}\label{th-Stric} Let $n\geq 2$ and $q\in (4,\infty)$ $(n=2)$, 
$q\in(2,\infty)$ $(n\geq 3)$. 
Let $\si:=n/2-1/q$. Then there exists $C>0$, which is independent of $T>0$, such that we have the inequality 
\begin{equation}\label{eq-2.1}
\|U(t)f\|_{L^q_T L^\infty_x}
\leq 
C\|f\|_{\Hs^{\si}({\mathbb R}^n)}.
\end{equation}
In particular, for any solutions to \eqref{eq-2.2}, we have
\begin{equation}\label{eq-2.3}
\|u(t)\|_{\Hs^s({\mathbb R}^n)}
\leq
\|u_0\|_{\Hs^s({\mathbb R}^n)}
+
\int_0^t
\|G(\tau)\|_{
\Hs^s({\mathbb R}^n)}d\tau,\,\,t>0, s\in [0, n/2)
\end{equation}
and,
\begin{equation}\label{eq-2.4}
\|u\|_{L^q_T L^\infty_x}
\leq
\|U(t)u_0\|_{L^q_TL^\infty_x}
+
C\|G\|_{L^1_T \Hs^{\si}},\,\,T>0.
\end{equation}
Furthermore, 
if $u_0(x)$ and $G(t,x)$ are radially symmetric about the origin $x=0$, 
then the estimates \eqref{eq-2.1} and \eqref{eq-2.4} remain true also for $q\in (2,4]$ $(n=2)$ 
and for $q=2$ $(n\ge 3)$. Here
and in what follows,
 we will use the notation
$\|f\|_{L^q_T X}:=\|f\|_{L^q(0,T; X)}$.
\end{lem}

\begin{prf} The Strichartz-type estimate \eqref{eq-2.1} 
is proved by 
Escobedo and Vega \cite{EsVe97} for $n=3$,
Klainerman and Machedon \cite{KM} for $n\ge 2$
(see also Proposition 1 of \cite{FW2}). 
The estimate \eqref{eq-2.3} is a consequence of 
the unitarity of $U(t)$ on $L^2({\mathbb R}^n)$. 
Also, the estimate \eqref{eq-2.4} follows directly from \eqref{eq-2.1}. 
In the presence of radial symmetry, 
the improvement of the range of admissible pairs $(q,r)$ 
was pointed out by Klainerman and Machedon 
for the solutions to the second order equation $\Box u=0$ 
in three space dimensions. 
They proved 
$$
\left\|D^{-1}\sin t D g\right\|_{L^2(0,\infty; L^\infty(\R^3))}
\leq
C\|g\|_{L^2({\mathbb R}^3)}
$$ 
for radially symmetric $g\in L^2$. 
In \cite{FW2}, Fang and Wang performed estimates 
directly for $U(t)f$, 
and proved \eqref{eq-2.1} for $q\in(2,\infty)$ $(n=2)$, 
$q\in [2,\infty)$ $(n\geq 3)$. 
\end{prf}

\begin{rem} As explained above, it is proved in \cite{FW2} that 
the inequality \eqref{eq-2.1} holds for $q=2$ and all $n\geq 3$ 
when $u_0$ is radially symmetric, 
which settles the conjecture made by Klainerman 
in \cite{Kla} (see Remark 1 there).
\end{rem}

\subsection{Local energy estimates }\label{sec-LE}
As is well-known, the local energy estimates (which are also known as KSS type estimates, Morawetz (radial) estimates) 
are indispensable for many problems for the wave equations. 
In this subsection, for solutions to the half-wave equations, we prove the analogs of 
the local energy estimates.
\begin{prp}\label{th-LE}
Let $n\geq 2$ and
$0<\delta\leq 1$.
There exists a constant $C>0$ 
depending only on $n$ and $\delta$, and independent of $T>0$,
such that any solutions 
to \eqref{eq-2.2} satisfy
\beeq\label{eq-3.1}
\|u\|_{L^\infty_T L^2_x}
+
T^{-\de/2}
\|
r^{-(1-\delta)/2}u
\|_{L^{2}_TL^2_x}
\leq
C\|u_0\|_{L^2({\mathbb R}^n)}
+
CT^{{\delta}/2}
\|
r^{(1-\delta)/2}G
\|_{L^{2}_T L^2_x},\eneq
\begin{align}
\label{eq-3.1'}
 & \|u\|_{L^\infty_T L^2_x}
+
(\ln(2+T))^{-1/2}
\|
r^{-(1-\delta)/2}\<r\>^{-\de/2}u
\|_{L^{2}_TL^2_x}
  \\
\leq    &  
C\|u_0\|_{L^2({\mathbb R}^n)}
+
C
(\ln(2+T))^{1/2}
\|
r^{(1-\delta)/2}\<r\>^{{\de}/2}G
\|_{L^{2}_T L^2_x}\ .\nonumber
\end{align}
In addition, 
if $0<\delta<\delta'$, $\delta\leq 1$, 
and $1+\delta'-\delta<n$,
there exists a constant $C>0$ depending on 
$n$, $\delta$, $\delta'$ such that 
the solution to  \eqref{eq-2.2} satisfies
\begin{align}
&
\|u\|_{L^\infty(0,\infty;L^2({\mathbb R}^n))}
+
\|
r^{-(1-\delta)/2}
\langle r\rangle^{-\delta'/2}
u
\|_{L^2((0,\infty)\times{\mathbb R}^n)}\\
&
\leq
C\|u_0\|_{L^2({\mathbb R}^n)}
+
C
\|
r^{(1-\delta)/2}
\langle r\rangle^{\delta'/2}
G
\|_{L^2((0,\infty)\times{\mathbb R}^n)}.\nonumber
\end{align}
\end{prp}

At first, let us record a version of local energy estimates for the standard wave operator $\Box=\pt^2-\Delta$. 
\begin{lem}\label{thm-LE}
Let $n\ge 2$, 
we have
\begin{equation}
\label{eq-LE1}
\|(\pt,\nabla_x) v\|_{L^\infty_t L^2_x \cap \dot \ell^{-1/2}_\infty L^2_t L^2_x }\les \|(\pt,\nabla_x) v(0)\|_{L^2}+\|\Box v\|_{L^1_t L^2+
 \dot \ell^{1/2}_1 L^2_t L^2_x 
}\ .
\end{equation} Here, with $q\in [1,\infty]$,
$
\|f\|_{\dot \ell^{s}_q L^2_t L^2_x}=\|2^{js} f\|_{\ell^q_{j\in\Z} L^2_{t,x}(t\in\R,|x|\simeq 2^j)}$.
In particular,
we have
\begin{align}
\label{eq-LE2}
&\|\pa v\|_{L^\infty_T L^2_x}
+
T^{-\delta/2}
\|
r^{-(1-\delta)/2}\pa v
\|_{L^{2}_TL^2_x}\\
\leq&
C\|\pa v(0)\|_{L^2({\mathbb R}^n)}
+
CT^{{\delta}/2}
\|
r^{(1-\delta)/2}\Box v
\|_{L^{2}_T L^2_x},
\nonumber
\end{align}
\begin{align}\label{eq-LE3}
 & \|\pa v\|_{L^\infty_T L^2_x}
+
(\ln(2+T))^{-{1}/2}
\|
r^{-(1-\delta)/2}\<r\>^{-{\de}/2}\pa v
\|_{L^{2}_TL^2_x}
  \\
\leq    &  
C\|\pa v(0)\|_{L^2({\mathbb R}^n)}
+
C
(\ln(2+T))^{1/2}
\|
r^{(1-\delta)/2}\<r\>^{\de/2}\Box v
\|_{L^{2}_T L^2_x}, \nonumber
\end{align}
\begin{align}\label{eq-LE4}
&
\|\pa v\|_{L^\infty(0,\infty;L^2({\mathbb R}^n))}
+
\|
r^{-(1-\delta)/2}
\langle r\rangle^{-\delta'/2}
\pa v
\|_{L^2((0,\infty)\times{\mathbb R}^n)}\\
&
\leq
C\|\pa v(0)\|_{L^2({\mathbb R}^n)}
+
C
\|
r^{(1-\delta)/2}
\langle r\rangle^{\delta'/2}
\Box v
\|_{L^2((0,\infty)\times{\mathbb R}^n)}.\nonumber
\end{align}
\end{lem}
The estimate \eqref{eq-LE1} for $n\ge 3$ could be proved by multiplier method,
see, e.g., \cite{MetSo06, W15}. The estimate for $n= 2$
is obtained in \cite{MeTa12MA}, even for small perturbation of the flat metric.
Local energy estimates have rich history and we refer \cite{MeTa12MA, LMSTW} for more exhaustive history of such estimates.

The estimates \eqref{eq-LE2}-\eqref{eq-LE4}
are also known as KSS type estimates \cite{KSS, HiYo06, MetSo06, Hi07, HWY1}, 
which are known to be implied by the local energy estimates \eqref{eq-LE1}, see, e.g., \cite{KSS, MetSo06, W15} and \cite[Lemma 1.5]{W15ext}.

To obtain \Cref{th-LE} from \Cref{thm-LE}, we need to check
that the weight functions, $w$, behave well under various operations, and 
it amounts to property $w^2\in A_2$.
\begin{lem}[Lemma 2.7 of \cite{HJLW}]\label{thm-weight}
Let $n\geq 1$ and let 
$w(x)
:=
r^{-(1-\delta)/2}
\< r\>^{-\delta'/2}$, 
with 
$0\leq 1-\delta\leq 1-\de+\delta'<n$. 
Then $w^2\in A_1$, where
 by $A_p$ we mean the Muckenhoupt $A_p$ class with $p\in [1,\infty)$. 
 In particular, it applies to all of the three weight functions occurred in \Cref{th-LE}. 
\end{lem}

\noindent {\bf Proof of Proposition \ref{th-LE}.} 
The idea is to reduce the required estimates to what are known for solutions to 
the second-order wave equation. 
Let $w_0\in A_1$ be any of the three weight functions occurring in Proposition \ref{th-LE}, and the estimates to be proven could be rewritten as follows
\begin{equation}\label{eq-LE-3.12}
\|
 u
\|_{L^{\infty}_TL^2_x}
+
\|
w_0 u
\|_{L^{2}_TL^2_x}
\leq   
C\|u_0\|_{L^2_x}+
C
\|
w_0^{-1} G
\|_{L^{2}_T L^2_x}, n\ge 2.
\end{equation}

Now,
for the solution $u$ 
to \eqref{eq-2.2}, 
let us introduce the auxiliary equation
\begin{equation}\label{eq-3.10}
i\partial_t v
+
\sqrt{-\Delta}v
=-u,
\quad
v(0)=0.
\end{equation}
Then, we see that $v$ satisfies
$$
(\partial_t^2-\Delta)v
=
-(i\partial_t-\sqrt{-\Delta})(i\partial_t+\sqrt{-\Delta})v
=
(i\partial_t-\sqrt{-\Delta})u=G
$$
and 
$$
\partial_t v(0)
=
-\frac{1}{i}
\bigl(
u(0)+\sqrt{-\Delta}v(0)
\bigr)=i u_0.
$$
Based on \Cref{thm-LE}, we get
\begin{equation}
\|
\pa v
\|_{L^{\infty}_TL^2_x}
+
\|
w_0 \pa v
\|_{L^{2}_TL^2_x}
\leq   
C\|u_0\|_{L^2_x}+
C
\|
w_0^{-1} G
\|_{L^{2}_T L^2_x}, n\ge 2.
\end{equation}
By \eqref{eq-3.10}, we see that the proof of \eqref{eq-LE-3.12} is reduced to 
\beeq\label{eq-LE-3.15}\|w_0 \sqrt{-\Delta}v \|_{L^2((0,T)\times{\mathbb R}^n)}
\les
\|w_0 \nabla v\|_{L^2((0,T)\times{\mathbb R}^n)}.
\eneq
However, as
$w_0^2\in A_1\subset A_2$ by \Cref{thm-weight},
we know that
$R_j=\sqrt{-\Delta}^{-1}\pa_j$,
the $j$-th Riesz transform, is bounded on $L^2(w_0^2 dx)$, and so
\begin{align}
&
\|
w_0\sqrt{-\Delta}w
\|_{L^2({\mathbb R}^n)}
=
\|
w_0
(-\Delta)^{-1/2}(-\Delta)
w
\|_{L^2({\mathbb R}^n)}\\
&
\leq
\sum_{j=1}^n
\|w_0
R_j\partial_j
w
\|_{L^2({\mathbb R}^n)}
\leq
C
\sum_{j=1}^n
\|w_0
\partial_j
w
\|_{L^2({\mathbb R}^n)},\nonumber
\end{align}
which completes the proof of \eqref{eq-LE-3.15}. The proof of  Proposition \ref{th-LE} is finished. 
{\hfill  {\vrule height6pt width6pt depth0pt}\medskip}

\section{Proof of Theorem 
\ref{th-1.1}}\label{sec-1.1}

The proof of 
\Cref{th-1.1} requires two ingredients. 
They are the $L_t^qL_x^\infty$ Strichartz-type estimates,
\Cref{th-Stric},
 and 
the chain rule of fractional orders,
\Cref{thm-Leib-h}.
We will give the proof for non-algebraic nonlinearity in the following, and the proof in the case of algebraic nonlinearity follows the  same lines of proof, 
with \eqref{eq-chain-3} replaced by \eqref{eq-Leib2}.

\subsection{Subcritically local wellposedness in $H^s$}
Define the nonlinear mapping
\begin{equation}\label{eq-4.1}
\Phi[u](t):=
U(t)u_0
-i
\int_0^t
U(t-\tau)
F(u(\tau))
d\tau.
\end{equation}
The existence of solution is equivalent to seeking a fixed point of $\Phi$. The proof of local wellposedness for $H^{s}$ with $s>n/2$ and $q=\infty$ is routine consequence of classical energy argument, and we will focus only on the case $s\in (0, n/2]$ with $q<\infty$ (and so is $n\ge 2$).

Thanks to the assumption $p>s$, $k=\lceil p\rceil-1$ and $q\ge p-1$, we get by \eqref{eq-2.3} and \eqref{eq-chain-4}, that, for any $t\in [0,T]$ and $s_{0}\in\{ 0, s\}$,
\begin{align}
\|\Phi[u](t)\|_{H^{s_{0}}}
&
\leq
\|u_0\|_{H^{s_{0}}}
+
C_{1}\int_0^t
\|u(\tau)\|_{L^\infty}^{p-1}
\|u(\tau)\|_{H^{s_{0}}}d\tau\label{eq-4.2}\\
&
\leq
\|u_0\|_{H^{s_{0}}}
+
C_{1}T^{1-\frac{p-1}q}\|u\|_{L^{q}_T L^\infty}^{p-1}
\|u\|_{L^\infty_T H^{s_{0}}}.\nonumber
\end{align}
As
 $q\in (\max({4}/(n-1), 2),\infty)$ and ${1}/{q}\ge n/2-s$, we have
 $H^{s}\subset \Hs^{n/2-1/q}$, and by \eqref{eq-2.1} and\eqref{eq-2.4},
\begin{eqnarray}
\|\Phi[u]\|_{L^{q}_T L^\infty_x}
&\leq&
\|U(t)u_0\|_{L^{q}_T L^\infty_x}
+
C_2T^{1-\frac{p-1}q}\|u\|_{L^{q}_T L^\infty_x}^{p-1}
\|u\|_{L^\infty_T H^s}\label{eq-4.3}\\
&\leq &
C_{3}\|u_0\|_{H^{s}}
+
C_2T^{1-\frac{p-1}q}\|u\|_{L^{q}_T L^\infty_x}^{p-1}
\|u\|_{L^\infty_T H^s}
.\label{eq-4.4}\end{eqnarray}
We also have
\beeq
\|\Phi[u](t)-\Phi[v](t)\|_{L^2}
\leq
C_4
T^{1-\frac{p-1}q}(
\|u\|_{L^{q}_T L^\infty_x}
+
\|v\|_{L^{q}_T L^\infty_x}
)^{p-1}
\|u(t)-v(t)\|_{L^\infty_T L^2}.\label{eq-4.5}\eneq

For the case $q>p-1$, setting $$T= c \|u_0\|_{H^{s}}^{-\frac{q(p-1)}{q-(p-1)}}
$$ for some small $c>0$ such that
$$C_{1}T^{1-\frac{p-1}q}(2C_{3}\|u_0\|_{H^{s}})^{p-1}\le 1/2,
C_2T^{1-\frac{p-1}q}(2C_{3}\|u_0\|_{H^{s}})^{p-1}\le C_{3}/2 \ ,$$
$$
C_4
T^{1-\frac{p-1}q}(4C_{3}\|u_0\|_{H^{s}})^{p-1}\le  1/2\ .$$
then we see that
$$\|\Phi[u]\|_{L^{\infty}_{T}H^{s}}
\leq
2\|u_0\|_{H^{s}},
\|\Phi[u]\|_{L^{q}_T L^\infty_x}\le 2C_{3}\|u_0\|_{H^{s}}\ ,$$
$$\|\Phi[u](t)-\Phi[v](t)\|_{L^2}
\leq
\frac{1}2\|u(t)-v(t)\|_{L^\infty_T L^2}.$$
That is, we have proved that
$\Phi$ is a contraction 
mapping on the complete metric space
\begin{align}
X_{1}:=\{u\in
L^\infty(0,T;H^{s})
\cap
C([0,T];L^2)
\cap 
L^{q}_{T}L^\infty_x,
u(0)=u_0,\\
\|u(t)\|_{H^{s}}
\leq
2\|u_0\|_{H^{s}},
\|u\|_{L^{q}_T L^\infty_x}\le 2C_{3}\|u_0\|_{H^{s}}
\}\nonumber
\end{align}
with the metric 
$d(u,v)
:=
\|u-v\|_{L^\infty_{T}L^2}$. The unique fixed point $u\in X_{1}$ is the unique solution we are looking for.

We can easily verify the fact that $u\in C([0, T]; H^s)$. Indeed, we know that $U(t)u_0\in C(\R; H^s)$
and
$\|F(u(\tau))\|_{H^s}\in L^1([0, T])$, which yields that 
$$\int_0^t
U(t-\tau)
F(u(\tau))
d\tau\in C([0,T]; H^s)\ ,$$
we find the right hand side of \eqref{eq-4.1} is in $C([0,T]; H^s)$ and so is
$u\in C([0, T]; H^s)$.

In addition, if $u_{0}\in H^{s_{1}}$,  
observing that
\eqref{eq-4.2} with $s_{0}=s_{1}$ also holds with some new constant $C_{1}'$. Then
see that 
$\Phi$ is also a contraction 
mapping on the complete metric space
\begin{align}
X_{2}:=\{u\in
L^\infty(0,T;H^{s_{1}})
\cap
C([0,T];L^2)
\cap 
L^{q}_{T}L^\infty_x,
u(0)=u_0,\\
\|u(t)\|_{H^{s_{1}}}
\leq
2\|u_0\|_{H^{s_{1}}},
\|u\|_{L^{q}_T L^\infty_x}\le 2C_{3}\|u_0\|_{H^{s}}
\}\nonumber
\end{align}
for the same kind of $T$ and
 the metric 
$d(u,v)
:=
\|u-v\|_{L^\infty_{T}L^2}$, which essentially
verifies the persistence of regularity property. The details are left to the interested reader.

\subsection{Critically local wellposedness in $H^{s_{c}}$}
For the critical case $s=s_{c}\in (0, n/2)$, we are forced to set $q=p-1\in (\max({4}/(n-1), 2),\infty)$. For this situation, 
\eqref{eq-4.2} and
\eqref{eq-4.3} 
 could be replaced by
\beeq
\|\Phi[u](t)\|_{\Hs^{s_{0}}}
\leq
\|u_0\|_{\Hs^{s_{0}}}
+
C_{1}\|u\|_{L^{p-1}_T L^\infty}^{p-1}
\|u\|_{L^\infty_T \Hs^{s_{0}}}\ ,
\forall t\in [0,T], s_{0}\in \{0, s_{c}\}\ ,
\label{eq-4.2'}
\eneq
\beeq
\|\Phi[u]\|_{L^{q}_T L^\infty_x}
\leq
\|U(t)u_0\|_{L^{q}_T L^\infty_x}
+
C_2\|u\|_{L^{p-1}_T L^\infty_x}^{p-1}
\|u\|_{L^\infty_T \Hs^s}\ .\label{eq-4.3'}\eneq

Let $A>0$ be a small constant satisfying
\begin{equation}
C_{1} (2A)^{p-1}\leq\frac12,\,\,
A+C_2(2A)^{p-1}(2\|u_0\|_{\Hs^{s_c}})
\leq
2A,\,\,
C_4(4A)^{p-1}\leq\frac12.
\end{equation}
By  the Strichartz estimate \eqref{eq-2.1}, 
$$\|U(t)u_0\|_{L^{p-1}(0,\infty; L^\infty_x)}
\leq C_{3}\|u_0\|_{\Hs^{s_c}}<\infty\ .$$
Thus, for this small constant $A$, 
we can choose $T_0>0$ so that 
\begin{equation}\label{eq-Th1small}
\|U(t)u_0\|_{L^{p-1}_{T_{0}} L^\infty_x}
\leq
A
\end{equation}
thanks to the absolute continuity of the integral. 
If we define
\begin{align}
X_{3}:=\{u=u(t,x)\,:\,
u\in
L^\infty_{T_0}H^{s_c}
\cap
C([0,T_0];L^2)
\cap 
L^{p-1}_{T_0}L^\infty_x,
u(0)=u_0,\\
\|u\|_{L^\infty_{T_0} L^2}\leq
2\|u_0\|_{L^2},
\| u\|_{L^\infty_{T_0} \Hs^{s_c}}\leq
2\|u_0\|_{\Hs^{s_c}},
\|u\|_{L^{p-1}_{T_0} L^\infty_x}
\leq 2A\}\nonumber
\end{align}
with the metric 
$d(u,v)
:=
\|u-v\|_{L^\infty_{T_0}L^2}$, 
then
 $X$ is a complete metric space and 
we see from  \eqref{eq-4.2'}, \eqref{eq-4.3'} and \eqref{eq-4.5}  that 
$\Phi$ is a contraction 
mapping on $X_{3}$. 

\subsection{Small data scattering in $H^{s_{c}}$}
When the initial data satisfies
$$\|u_0\|_{\Hs^{s_c}}\le A/C_{3}\ ,$$
we see that we could set $T_{0}=\infty$ and so is the
global existence of $u\in 
L_{t}^\infty H^{s_c}
\cap
C_{t}L^2
\cap 
L^{p-1}_{t}L^\infty_x$. It remains to prove scattering for the global solution.

By \Cref{thm-Leib-h}, $F(u)\in
L^1(0,\infty; H^{s_c})$.
Using the basic inequality for vector-valued integrable functions, 
we have
$$\|\int_{t}^{\infty}
U(-\tau)F(u(\tau))d\tau\|_{H^{s_{c}}}\le C
\|F(u)\|_{L^1(t,\infty; H^{s_{c}})},
$$ which converges to $0$ as $t\to +\infty$.
In other words, we have proven that
$$u_0
-i
\int_0^t
U(-\tau)
F(u(\tau))
d\tau
\to
u_0-i
\int_{0}^{\infty} U(-\tau)F(u(\tau))d\tau:=u_{0}^{+}\ ,
$$
in $H^{s_{c}}$.
Thus,
we see that
$$\|u(t)-U(t)u_0^+\|_{H^{s_c}}=
\|U(t)\int_{t}^{\infty} U(-\tau)F(u(\tau))d\tau\|_{H^{s_c}}=\|\int_{t}^{\infty} U(-\tau)F(u(\tau))d\tau\|_{H^{s_c}}\to 0
,$$
as ${t\to +\infty}$. This completes the proof of \Cref{th-1.1}.


\section{Proof of Theorems \ref{th-1.2} and \ref{th-1.3}}\label{sec-1.2}
In this section, we prove Theorems \ref{th-1.2} and \ref{th-1.3}, as well as the existence part in  \Cref{th-1.4}.
The key to the proof is the local energy estimates \Cref{th-LE}, as well as the radial Sobolev inequalities, and the following weighted fractional chain rule of  \cite[Theorem 2.5]{HJLW}
\begin{lem}[Fractional chain rule]\label{thm-Leibniz}
Let $s\in (0,1)$, $q\in (1,\infty)$. Assume $F:\R^k\rightarrow \R^l$ is a $C^1$ map, satisfying $F(0)=0$ and 
\beeq\label{eq-wLeib-assu}|F'(\tau v+(1-\tau)w)|\le \mu(\tau)|G(v)+G(w)|,\eneq
with $G>0$ and $\mu\in L^1([0,1])$.
If
$w_1^{q}, (w_{1}w_2)^{q}\in A_q$ and $w_2^{-1}\in A_1$, we have
\beeq\label{eq-wLeib}\|w_1 w_{2} D^s F(u)\|_{L^q}\les \|w_1 D^s u\|_{L^q} \|w_2 G(u)\|_{L^\infty}\ .\eneq
In particular, for $s\in [0,1]$,
if $w^{2}\in A_1$, then we have
\beeq\label{eq-wLeib2}\|w^{-1} D^s F(u)\|_{L^2}\les \|w D^s u\|_{L^2_x} \|w^{-2} G(u)\|_{L^\infty_x}\ .\eneq
\end{lem}
We notice that we have added 
the trivial cases $s=0, 1$ to \eqref{eq-wLeib2}, due to the fact that
if $w^{2}\in A_{2}$, $$\|w D u\|\simeq \|w \nabla u\|\ .$$



\subsection{Proof of Theorem \ref{th-1.2}}
Let
$n\geq 2$, $p\in (1,1+2/(n-1))$, $\de=1-\frac{n-1}{2}(p-1)\in (0,1)$ and $w(x):=|x|^{-(1-\de)/2}$.
By \Cref{thm-weight}, $w^{2}\in A_{1}$ and we can apply 
\eqref{eq-wLeib2} with $G(u)=|u|^{p-1}$: for $s\in [0,1]$,
$$\|w^{-1} D^s F(u)\|_{L^2}\les \|w D^s u\|_{L^2_x} \|w^{-2} G(u)\|_{L^\infty_x}
\les \|w D^s u\|_{L^2_x} \|w^{-2/(p-1)} u\|_{L^\infty_x}^{p-1}
\ .$$
Then,
applying  
 \eqref{eq-3.1}, we see that the map  $\Phi$ defined in \eqref{eq-4.1} satisfies
$$\|\Phi[u]\|_{X^{s}}:=\|\Phi[u]\|_{L^\infty_T \Hs^{s}_x}
+
T^{-{\delta}/2}
\|
w D^{s} \Phi[u]
\|_{L^{2}_TL^2_x}
\les
\|u_0\|_{\Hs^{s}({\mathbb R}^n)}
+
T^{{\delta}/2}
\|
w^{-1}D^{s}G
\|_{L^{2}_T L^2_x}$$
and so, for
fixed $s\in (1/2,1]\cap (1/2,n/2)$, $s_{1}\in (1/2, s]$, we have
for $\sigma=0,1-s_1,s_1,s$,
\beeq\label{eq-5.4}
\|\Phi[u]\|_{X^{\si}}
\le C_{5}(
\|u_0\|_{\Hs^{\si}({\mathbb R}^n)}
+
T^{\frac{\delta}2}
\|
w D^{\si} u
\|_{L^{2}_TL^2_x}
\|w^{-\frac{2}{p-1}} u\|_{L^\infty_T L^\infty_x}^{p-1})\ .
\eneq

Note that 
$w^{-{2}/(p-1)}=
|x|^{(1-\delta)/(p-1)}=|x|^{(n-1)/2}$. 
Therefore we can use the radial Sobolev lemma of \cite{HJLW} 
(see (2.2) of \cite{HJLW} 
which is a direct consequence 
of the ``end-point'' trace lemma of \cite{FW11} 
and the real interpolation) 
to get
\begin{equation}\label{eq-5.5}
\|
r^{(n-1)/2}u
\|_{L^\infty({\mathbb R}^n)}
\leq
C
\|
D^{s_1} u
\|_{L^2({\mathbb R}^n)}^{1/2}
\|
D^{1-s_1} u
\|_{L^2({\mathbb R}^n)}^{1/2}.
\end{equation}
It is at this place that we need the assumption of 
radial symmetry in our proof.

Let us define the space of functions 
where we carry out the contraction-mapping argument: 
we set 
\begin{align}
X(s,s_1,T):=
\{
u=u(t,x),\,\mbox{{\rm radially symmetric about $x=0$}} :\\
\|
u
\|_{X^{\si}}
\leq
2C_5\|u_0\|_{\Hs^{\sigma }({\mathbb R}^n)}
\,\,\mbox{for}\,\sigma=0,1-s_1,s_1,s\}.\nonumber
\end{align}
Here, $C_5$ is the constant appearing in \eqref{eq-5.4}. 
Equipped with the metric 
$d(u,v)
:=T^{-\delta/2}
\|
w(u-v)
\|_{L^2((0,T)\times{\mathbb R}^n)}
+
\|
u-v
\|_{L^\infty(0,T;L^2({\mathbb R}^n))}$, 
$X(s,s_1,T)$ is a complete metric space. 
By \eqref{eq-5.4} and \eqref{eq-5.5}, we find that 
the mapping $\Phi[u]$ satisfies for 
$\sigma=0,1-s_1,s_1,s$
\beeq
\|
\Phi[u]
\|_{X^{\si}}
\leq
C_5\|u_0\|_{\Hs^\si({\mathbb R}^n)}
+
C T^\delta
\|u\|_{X^{\si}}
(\|u\|_{X^{s_{1}}}\|u\|_{X^{1-s_{1}}})^{(p-1)/2}\ ,\eneq
\beeq
\|
\Phi[u]-\Phi[v]
\|_{X^{0}}
\leq
C T^\delta
\|u-v\|_{X^{0}}
(\|u\|_{X^{s_{1}}}\|u\|_{X^{1-s_{1}}}+\|v\|_{X^{s_{1}}}\|v\|_{X^{1-s_{1}}})^{(p-1)/2}\eneq
for any $u, v\in X(s,s_1,T)$. 
Choosing $T>0$ as in \eqref{eq-life-lower} for some small constant $c>0$,  then
$$T^{\de}
(\|D^{s_1}
u_0
\|_{L^2({\mathbb R}^n)}^{1/2}
\|
D^{1-s_1}
u_0
\|_{L^2({\mathbb R}^n)}^{1/2})^{p-1}\ll 1$$
and
we easily see that 
$\Phi[u]$ maps $X(s,s_1,T)$ into itself 
and it is a contraction-mapping there. 
In this way, we can prove the local wellposedness 
in $H_{\mbox{\scriptsize rad}}^s$ when $n\geq 2$, 
$s\in (1/2,1]\cap (1/2,n/2)$ and $0<p-1<2/(n-1)$. 
The proof of   \Cref{th-1.2} is finished. 
\subsection{Small data long time existence, Theorems \ref{th-1.3} and \ref{th-1.4}}
The proof of small data global existence, Theorem \ref{th-1.3}, as well as the 
small data long time existence part of Theorem \ref{th-1.4}, proceeds along the same lines, with trivial modifications.

More precisely, 
let
$$w:= \left\{\begin{array}{ll }
     r^{-(1-\delta)/2}\<r\>^{- \de/2}&    p=1+2/(n-1)\\
     r^{-(1-\delta)/2}\<r\>^{-\de'/2}&    p\in(1+2/(n-1), 1+2/(n-2))
\end{array}\right.$$
and
$$A_{T}:= \left\{\begin{array}{ll }
    ( \ln(2+T))^{1/2}&    p=1+2/(n-1)\\
 1&    p\in(1+2/(n-1), 1+2/(n-2))\ .
\end{array}\right.$$
Here,
we choose $\de'>\de$ to be fixed such that \Cref{thm-weight} applies for $w$.
Then we have for any $T\in (0,\infty)$ and $\si\in [0,1]$,
\begin{eqnarray*}
\|\Phi[u]\|_{X^{\si}}&
\les&
\|u_0\|_{\Hs^{\si}({\mathbb R}^n)}
+
A_{T}
\|
w D^{\si} u
\|_{L^{2}_TL^2_x}
\|w^{-\frac{2}{p-1}} u\|_{L^\infty_T L^\infty_x}^{p-1}\\
&\les&
\|u_0\|_{\Hs^{\si}({\mathbb R}^n)}
+
A_{T}^{2}
\|
u
\|_{X^{\si}}
\|w^{-\frac{2}{p-1}} u\|_{L^\infty_T L^\infty_x}^{p-1}
\end{eqnarray*}
and
$$
\|\Phi[u]-\Phi[v]\|_{X^{0}}
\les
A_{T}^{2}
\|
u-v
\|_{X^{0}}
\|w^{-\frac{2}{p-1}} (u,v)\|_{L^\infty_T L^\infty_x}^{p-1}
$$
where
 $$\|u\|_{X^{s}}:=\|u\|_{L^\infty_T \Hs^{s}_x}
+
A_{T}^{-1}
\|
w D^{s} u
\|_{L^{2}_TL^2_x}\ .$$
Based on these estimates, the proof of 
small data long time existence is easily reduced to the 
 radial Sobolev lemma, \eqref{eq-5.5} and
\begin{equation}\label{eq-5.5'}
\|r^{n/2-s}  u\|_{L^\infty_{x}}\les \| u\|_{ \Hs^s}, s\in (1/2, n/2)
\end{equation}
(see, e.g.,  \cite{FW11}, \cite[(2.2)]{HJLW}).

\section{Nonexistence of global solutions and upper bound of the lifespan}\label{sec-bu}
In this section, we give the remaining part of  Theorem \ref{th-1.4},
nonexistence of global solutions and upper bound of the lifespan, for the case of $1<p\le 1+2/(n-1)$ and $F(u)=i |\Re u|^p$, that is
\beeq\label{eq-6.1}i\partial_t u-\sqrt{-\Delta}u=i |\Re u|^p\ .\eneq

As we have mentioned in the introduction, the idea is to establish a connection between the half-wave problems and the nonlinear wave equations. The procedure is similar to that of 
\Cref{th-LE}.
Introduce $v$ in the way similar to \eqref{eq-3.10}, that is,
\beeq\label{eq-6.2}i\partial_t v
+
\sqrt{-\Delta}v
=i u,
\quad
v(0)=0.
\eneq
Then, we see that $v$ satisfies
\beeq
(\partial_t^2-\Delta)v
=
-(i\partial_t-\sqrt{-\Delta})(i\partial_t+\sqrt{-\Delta})v
=-i
(i\partial_t-\sqrt{-\Delta})u= |\Re u|^p
\eneq
and 
$$v(0)=0,
\partial_t v(0)
=
\frac{1}{i}
\bigl(
iu(0)-\sqrt{-\Delta}v(0)
\bigr)= u_0.
$$

For any real-valued, radial function $g\in C_0^\infty(\R^n)$ with $\int gdx>0$, 
we set the initial data to be $u_{0}=\ep g$ with $\ep>0$.
Then we see that $v$ is  real-valued and so is $\sqrt{-\Delta} v$ (as $\overline{\sqrt{-\Delta}v}=
\sqrt{-\Delta}\bar v$).
Thus
\beeq\label{eq-6.4.1}
-i\partial_t v
+
\sqrt{-\Delta}v
=-i \bar u, \ \pt v=\frac{u+\bar u}2=\Re u, \sqrt{-\Delta} v=i\frac{u-\bar u}2=-\Im u
\ ,\eneq
and so
\beeq\label{eq-6.4}
(\partial_t^2-\Delta)v=|\pt v|^{p}, v(0)=0,
\partial_t v(0)
=
\ep g\ .\eneq

In summary, we have seen that for any solution $u$
to \eqref{eq-6.1} with data $u_{0}=\ep g$, the solution $v$  to \eqref{eq-6.2} satisfies
\eqref{eq-6.4}.
Then it is well-known that, for $p\in (1, 1+2/(n-1)]$,
\eqref{eq-6.4} could not admit global solutions for any $\ep>0$, and the 
lifespan $T$ satisfies $T\le L_{\ep}$ with $L_{\ep}$ given by \eqref{eq-upperbd}, for some constant $C>0$,
see, e.g., \cite{Zhou01} and references therein.

For completeness, we give here a simple proof by using the method of test functions, in the current setting, as 
appeared in \cite{YorZh06}, \cite{ZhouHan11} for wave equations.

Assume that $u\in
L^\infty(0,T;H^s)\cap 
C([0,T];L^2)\subset
C([0,T];H^{s-\de})$ (for any $\de>0$) is a solution to the equation 
\eqref{eq-6.1} with data $u_{0}=\ep g$. Then our task is to show that $T< L_{\ep}$.
As $v$ is the solution   to \eqref{eq-6.2}, we have
$v\in L^\infty(0,T;H^s)\cap C([0,T];H^{s-\de})$. Recalling \eqref{eq-6.4.1}, we have
\beeq
v\in 
L^\infty(0,T;H^{s+1})\cap 
C([0,T];H^{s+1-\de})\cap C^{1}([0,T];H^{s-\de})\ .\eneq

Recall  that, by using the similar proof of \Cref{th-1.1}
with radial Strichartz estimates \Cref{th-Stric}, the problem \eqref{eq-6.4} with radial data is locally wellposed in $L^\infty H^{s+1}\cap Lip H^s$ with $s\in [s_c, p)$ and $s\ge (n-1)/2$ (except $s=1/2$ and $n=2$). As we have $p\le 1+2/(n-1)$, $s_c\le 1/2$ and so the problem is locally wellposed  for
$s\in (1/2, p)$ for $n=2$ and 
$s\in [(n-1)/2, p)$ for $n=3, 4$.
In particular, 
there is $R>0$ such that $g$ is supported in $B_R$, and so 
$$\mathrm{supp}\ v\subset \{
(t,x)\in [0,T]\times\R^n, |x|\le t+R
\}\ .$$

Let
$$\phi(t,x)=\int_{\Sp^{n-1}} e^{x\cdot \om-t} d\om\ .
$$  Then,
it satisfies
$$\Delta\phi=\phi,
\pt \phi=-\phi,
\pt^2 \phi=\phi\ .$$
Moreover, it is well-known that
\beeq\label{eq-simeq} 0<
c (1+|x|)^{-\frac{n-1}2}e^{|x|-t}
\le \phi\le C (1+|x|)^{-\frac{n-1}2}e^{|x|-t} \ .\eneq
For the convenience of the reader, let us give a proof. The estimate for $r:=|x|\les 1$ is trivial and we need only to consider $r\gg 1$. In this case,
\begin{eqnarray*}
r^{\frac{n-1}2}e^{t-r}
\phi(t,x) & = & |\Sp^{n-2}|
\int_{0}^\pi e^{r (\cos \theta-1)} r^{\frac{n-1}2}(\sin\theta)^{n-2} d\theta\\
 & = & |\Sp^{n-2}|
\int_{0}^{\pi/2} e^{r (\cos \theta-1)} r^{\frac{n-1}2}(\sin\theta)^{n-2} d\theta+o(1)\ .
\end{eqnarray*}
  Then, as $2\theta/\pi\le \sin\theta\le \theta$ 
  and $ \theta^2/\pi\le
  1-\cos \theta\le \theta^2/2$ 
  for $\theta\in [0, \pi/2]$, we have
$$  \int_{0}^{\frac \pi 2} e^{r (\cos \theta-1)} r^{\frac{n-1}2}(\sin\theta)^{n-2} d\theta\le
 \int_{0}^{\pi/2} e^{-r \frac{\theta^2}{\pi}} r^{\frac{n-1}2}\theta^{n-2} d\theta 
\le \left(\frac{\pi}{2}\right)^{n-1}  \int_{0}^{\infty} e^{- \frac{\theta^2}{2}} \theta^{n-2} d\theta \ ,
$$
$$  \int_{0}^{\frac \pi 2} e^{r (\cos \theta-1)} r^{\frac{n-1}2}(\sin\theta)^{n-2} d\theta\ge
 \int_{0}^{\frac \pi 2} e^{-r \frac{\theta^2}{2}} r^{\frac{n-1}2}\left(\frac{2\theta}{\pi}\right)^{n-2} d\theta 
=
\left(\frac{2}{\pi}\right)^{n-2}   \int_{0}^{\frac{\pi\sqrt{r}}2} e^{- \frac{\theta^2}{2}}\theta^{n-2} d\theta \ ,
$$
  and so is \eqref{eq-simeq}.

Introducing a cutoff function $\psi$ which is identity on $B_{T+R}$, then we have $v_t=\psi v_t\in C([0,T]; H^{s-\de})\subset C([0,T]; L^p)$ for $p\in (1, 1+2/(n-1)]$,
$\int \phi |v_t|^p dx
=\int \phi |\psi v_t|^p dx
\in C([0,T])$.
As $\Delta v\in  C([0,T]; H^{s-1-\de})$, 
by \eqref{eq-6.4}, we have
$$\int v_{tt}\phi dx
=\int \phi \psi(\Delta v +|v_t|^p) dx\in C([0,T])
 $$
 and so
 $\int \phi v dx\in C^2$,
  $\int \phi v_t dx\in C^1$.

With $\Box=\pt^2-\Delta$,  we have
$\Box \phi=0$, 
$$
\phi \Box v -v\Box \phi=
\frac{d}{dt} (\phi\pt v-v\pt \phi)-\nabla\cdot (\phi\nabla v-v\nabla \phi)
,$$
and so
\beeq
F'=\frac{d}{dt} \int \phi(\pt v+v)dx=
\frac{d}{dt} \int (\phi\pt v-v\pt \phi)dx
=\int \phi \Box v -v\Box \phi dx=
\int\phi|v_t|^pdx
\ ,
\eneq
where $F(t):=\int \phi(\pt v+v)dx\in C^1$.
On the other hand,
let
$H(t):=\int \phi \pt v dx\in C^1$, 
\beeq
\int\phi \Box v -v\Box \phi dx
=\int(\phi  v_{tt} -v\phi_{tt} )dx=
H'-
\int (\phi_t  v_{t}+v\phi) dx
=H'+\int \phi  (v_{t}- v)dx
\ .
\eneq
Let
$G:=\int \phi (v_t-v)dx=2H-F$. Then we have
$F'=H'+G\ge 0$, 
$$ G'=2H'-F'\ge 2H'-2F'=-2G
$$
and so
$$G(t)\ge G(0) e^{-2t}=\ep e^{-2t} \int g\phi d x> 0\ .$$
Thus, as $F'(t)=\int \phi|v_t|^p dx$, 
$$F\le 2H\les
\left(\int \phi|v_t|^p dx\right)^{1/p}
\left(\int_{B_{t+R}} \phi dx\right)^{1-1/p}\les
F'(t)^{1/p} (t+R)^{\frac{n-1}2(1-1/p)},
$$
that is, there exists $c>0$ such that
\beeq\label{eq-6.10} F'(t)\ge c F^p (t+R)^{-\frac{(n-1)(p-1)}2}, \forall t\in [0, T],
F(0)=\ep \int g\phi\ge c\ep
\ .\eneq
Based on \eqref{eq-6.10}, it is a standard argument to conclude $T< L_\ep$ for some $C>0$.
The proof is finished.

\subsection*{Acknowledgment}
 The second author, C. W., is grateful to Professor
 Yoshio Tsutsumi for helpful discussion and valuable comments on the 
fractional chain rule.
K. H. was supported in part by JSPS KAKENHI Grant Number JP15K04955.
C. W. was supported in part by National Support Program for Young Top-Notch Talents.

\end{document}